\documentclass[11pt]{amsart}

\usepackage[a4paper,hmargin=3.5cm,vmargin=4cm]{geometry}
\usepackage{amsfonts,amssymb,amscd,amstext,verbatim}

\usepackage{fancyhdr}
\usepackage[utf8]{inputenc}
\usepackage{hyperref}
\usepackage{verbatim}

\usepackage{times}

\usepackage{enumerate}
\usepackage{titlesec}
\usepackage{mathrsfs}

\titleformat{\section}
{\filcenter\bfseries\large} {\thesection{.}}{0.2cm}{}
\titleformat{\subsection}[runin]
{\bfseries} {\thesubsection{.}}{0.15cm}{}[.]
\titleformat{\subsubsection}[runin]
{\em}{\thesubsubsection{.}}{0.15cm}{}[.]

\usepackage[up,bf]{caption}

\newtheorem{theorem}{Theorem}[section]

\newtheorem{lemma}[theorem]{Lemma}
\newtheorem{corollary}[theorem]{Corollary}

\theoremstyle{definition}

\newtheorem{remark}[theorem]{Remark}

\newtheorem{problem}[theorem]{Problem}

\numberwithin{equation}{section}
\numberwithin{figure}{section}

\newcommand\Ocal{\mathcal{O}}

\newcommand\Vcal{\mathcal{V}}
\newcommand\Wcal{\mathcal{W}}

\newcommand\C{\mathbb{C}}

\newcommand\N{\mathbb{N}}

\newcommand\igot{\mathfrak{i}}

\renewcommand\igot{\mathfrak{i}}

\renewcommand\imath{\igot}

\newcommand\hra{\hookrightarrow}

\newcommand\longhookrightarrow{\ensuremath{\lhook\joinrel\relbar\joinrel\rightarrow}}

\newcommand\wt{\widetilde}
\newcommand\wh{\widehat}

\newcommand\Aut{\mathrm{Aut}}

\newcommand{\bea}{\begin{eqnarray*}}
\newcommand{\eea}{\end{eqnarray*}}

\numberwithin{equation}{section}

\begin{document}

\fancyhead[LO]{Runge tubes in Stein manifolds with the density property}
\fancyhead[RE]{F.\ Forstneri{\v c} and E.\ F.\ Wold} 
\fancyhead[RO,LE]{\thepage}

\thispagestyle{empty}

\vspace*{1cm}
\begin{center}
{\bf\LARGE Runge tubes in Stein manifolds with \\  the density property}

\vspace*{0.5cm}

{\large\bf  Franc Forstneri{\v c} and Erlend F.\ Wold} 
\end{center}

\vspace*{1cm}

\begin{quote}
{\small
\noindent {\bf Abstract}\hspace*{0.1cm}
In this paper we give a simple proof of the existence and plenitude of Runge tubes in $\C^n$ $(n>1)$
and, more generally, in Stein manifolds with the density property. We show in particular that
for any algebraic submanifold $X$ of codimension at least two in a complex Euclidean space 
$\C^n$, the normal bundle of $X$ admits a holomorphic embedding onto a Runge domain in $\C^n$ 
which agrees with the inclusion map $X\hra\C^n$ on the zero section.

\vspace*{0.2cm}

\noindent{\bf Keywords}\hspace*{0.1cm}  Runge tube, holomorphic automorphism, Stein manifold, 
affine algebraic manifold, density property

\vspace*{0.1cm}

\noindent{\bf MSC (2010):}\hspace*{0.1cm}}  32E10; 32E30; 32H02; 32M17; 14R10

\end{quote}

\section{Introduction}

It has been an open question for a long time whether it is possible to embed $\C^*\times\C$ 
as a Runge domain in $\C^2$. (Here, $\C$ is the complex plane and $\C^*=\C\setminus\{0\}$.)
Such hypothetical domains have been called {\em Runge cylinders} in $\C^2$.
The question arose in connection with the classification of Fatou-components for 
H\'{e}non maps by E.\ Bedford and J.\ Smillie in 1991,  \cite{BedfordSmillie1991}. 
This problem has recently been solved in the affirmative 
by F.\ Bracci, J.\ Raissy and B.\ Stens\o nes \cite{BracciRaissyStensones2017X} 
who obtained a Runge embedding of $\mathbb C^*\times\mathbb C$ into $\mathbb C^2$
as the basin of attraction of a (non-polynomial) holomorphic automorphism of $\C^2$
at a parabolic fixed point.

The purpose of this note is to give a very simple proof of the existence of Runge cylinders, 
and furthermore of the existence of an abundance of Runge tube domains in all Stein manifolds with 
the density property.  Although our proof is completely different from that in \cite{BracciRaissyStensones2017X}, 
both proofs depend crucially on the Anders\'{e}n-Lempert theory (see \cite[Chapter 4]{Forstneric2017E}).  

Recall (see D.\ Varolin \cite{Varolin2000,Varolin2001} or \cite[Definition 4.10.1]{Forstneric2017E})
that a complex manifold $Y$ has the {\em density property} if every holomorphic vector field
on $Y$ is a  uniform limit on compacts of finite sums of $\C$-complete holomorphic 
vector fields on $Y$. In particular, Euclidean spaces $\C^n$ of dimension $n>1$ have the density
propery by E.\ Anders{\'e}n and L.\ Lempert \cite{AndersenLempert1992}.

The following is our first main result. 

%
%
\begin{theorem}\label{th:main}
Let $X$ and $Y$ be Stein manifolds with $\dim X<\dim Y$, and assume that $Y$ has
the density property.  Suppose that $\theta:X\hra Y$ is a holomorphic embedding with $\Ocal(Y)$-convex image 
(this holds in particular if $\theta$ is proper), and let $E\to X$ denote the normal bundle associated to
$\theta$. Then, $\theta$ is approximable uniformly on compacts in $X$ by 
holomorphic embeddings of $E$ into $Y$ whose images are Runge domains.  
\end{theorem}

Recall that a locally closed subset $Z$ of a complex manifold $Y$ is said to be {\em $\Ocal(Y)$-convex}
if for every compact set $K\subset Z$, its $\Ocal(Y)$-convex hull 
\begin{equation}\label{eq:hull}
	\wh K_{\Ocal(Y)}=\{y\in Y: |f(y)|\le \sup_K |f|\ \ \forall f\in\Ocal(Y)\}
\end{equation}
is compact  and contained in $Z$.

To get a Runge embedding of $\mathbb C^*\times\mathbb C$ into $\mathbb C^2$ from 
Theorem \ref{th:main}, one embeds $X=\mathbb C^*$  onto the curve $\{zw=1\}\subset \mathbb C^2$
and notes that any vector bundle over $\mathbb C^*$ (and in fact over any open Riemann surface) is trivial
by Oka's theorem \cite{Oka1939}. (See also \cite[Sect.\ 5.2]{Forstneric2017E}.)

There is a big list of Stein manifolds, and in particular of affine
algebraic manifolds, which are known to have the density property; see 
the list of examples in \cite{AndristForstnericRitterWold2016}, as well as 
the recent surveys of Kaliman and Kutzschebauch \cite{KalimanKutzschebauch2015} and of the first named
author \cite[Sect.\ 4.10]{Forstneric2017E}. The authors, together with R.\ Andrist and T.\ Ritter, proved in 
\cite{AndristForstnericRitterWold2016,AndristWold2014} that every Stein manifold 
$X$ embeds properly holomorphically into any Stein manifold $Y$ with the
density property satisfying $\dim Y\geq 2\dim X +1$. Every open Riemann surface $X$
embeds properly holomorphically into $\C^3$, and a plenitude of them
embed properly into $\C^2$; see \cite{ForstnericWold2013} and \cite[Sects.\ 9.10-9.11]{Forstneric2017E} 
for a discussion of this topic. By Theorem \ref{th:main}, every such embedding can be approximated by a 
Runge embedding of the normal bundle of $X$ into $Y$. This provides a huge variety of nontrivial
Runge tubes in any Stein manifold with the density property. 
In particular, we have the following corollary to Theorem \ref{th:main}.

%
%
\begin{corollary}[Runge tubes over open Riemann surfaces] \label{cor:tubesC2} 
If $X$ is an open Riemann surface which admits a proper holomorphic embedding into $\C^2$,
then $X\times \C$ admits a Runge embedding into $\C^2$. For every open Riemann surface $X$
and every $k\ge 2$, $X\times\C^k$  admits a Runge embedding into $\C^{k+1}$, and 
into any Stein manifold $Y^{k+1}$ with the density property.
\end{corollary}

%
%
\begin{remark}
The Runge embeddings $E\hra Y$ of the normal bundle in Theorem \ref{th:main} need not agree with 
the given embedding $\theta:X\hra Y$ on the zero section of $E$ (which we identify with $X$). 
Indeed, it is known that for every pair of integers $1\le k<n$ there exists a proper holomorphic
embedding $\theta:X=\C^k\hra Y=\C^n$ whose complement is $(n-k)$-hyperbolic 
in the sense of Eisenman; in particular, there are no nondegenerate holomorphic maps
$\C^{n-k}\to \C^n\setminus \theta(\C^k)$
(see  Buzzard and Forn{\ae}ss \cite{BuzzardFornaess1996} for the case $k=1,n=2$ 
and Borell and Kutzschebauch\cite{BorellKutzschebauch2006} and 
Forstneri\v c \cite[Corollary 5.3]{Forstneric1999} for the general case).
Since the normal bundle of the embedding $\theta$ is the trivial bundle $E=\C^n=\C^k\times \C^{n-k}\to \C^k$ 
and the complement of  $\C^k\times \{0\}^{n-k}$ in $\C^n$ is clearly
not  $(n-k)$-hyperbolic, $\theta$ does not extend to a holomorphic embedding $E=\C^n\hra \C^n$.
\qed\end{remark}

However, we can ensure this additional interpolation condition for algebraic embeddings 
of codimension at least 2 into $\C^n$. Here is the precise result.

%
%
\begin{theorem}\label{th:main2}
Let $X$ be a Stein manifold and $\theta\colon X\hra \C^n$ be a proper holomorphic embedding onto an
algebraic submanifold $\theta(X)\subset \C^n$. If $n\ge \dim X+2$, then $\theta$ extends to a 
holomorphic Runge embedding $\tilde\theta\colon E\hra Y$ of the normal bundle of $\theta$. 
\end{theorem}

According to Docquier and Grauert \cite{DocquierGrauert1960} 
(see also \cite[Theorem 3.3.3]{Forstneric2017E}), every proper holomorphic embedding $\theta\colon X\hra Y$ 
of a Stein manifold $X$ into a complex manifold $Y$ extends to an embedding of a neighborhood 
of the zero section in the normal bundle of $\theta$.  
Theorem \ref{th:main2} says that one can in fact embed the entire normal bundle 
as a Runge domain in $\C^n$.

The proof of Theorem \ref{th:main2} is similar to that of Theorem \ref{th:main}. 
It uses the result of Kaliman and Kutzschebauch \cite[Theorem 6]{KalimanKutzschebauch2008IM} 
that the Lie algebra of all algebraic vector fields on $\C^n$
vanishing on an algebraic submanifold of codimension at least two 
enjoys the algebraic density property. It follows that flows of such vector fields 
can be approximated by automorphisms of $\C^n$ fixing the submanifold pointwise. 

%
%
\begin{corollary}\label{cor:algebraic}
If $X$ is an affine algebraic curve, then every proper algebraic embedding $\theta\colon X\hra\C^{n+1}$
$(n\ge 2)$ extends to a holomorphic embedding $\tilde\theta:X\times \C^{n}\hra \C^{n+1}$
onto a Runge domain in $\C^{n+1}$. 
\end{corollary}

\begin{remark}
Note that holomorphic Runge embeddings of the normal bundle, furnished by
Theorem \ref{th:main2} and Corollary \ref{cor:algebraic}, can never be algebraic.
Indeed, if $E\to X$ is the algebraic normal bundle of an algebraic submanifold $X\subset\C^n$
and $F\colon E\hra \C^n$ is an algebraic embedding, then $\Omega=F(E)\subset \C^n$
is a Zariski open set in $\C^n$ and its complement $A=\C^n\setminus \Omega$ is a Zariski 
closed set, i.e., an algebraic subvariety of $\C^n$ (see Chevalley \cite{Chevalley1958}).
Since $\Omega$ is a Stein domain, $A$ must be of pure codimension one, and hence $A=\{f=0\}$
for some entire function $f\in \Ocal(\C^n)$. Clearly, the function $1/f\in \Ocal(\Omega)$ cannot 
be approximated uniformly on compacts in $\Omega$ by entire functions, and hence the domain $\Omega=F(E)$ 
is not Runge in $\C^n$.
\end{remark}

We conclude this introduction by pointing out the following open problem related to 
Theorem \ref{th:main} and Corollary \ref{cor:algebraic}.

\begin{problem}
Is there a Runge embedding of the (trivial) normal bundle $E=H\times \C\cong \C^*\times \C$
of the hyperbola $H = \{(z,w) \in \C^2: zw=1\}$ extending the inclusion map $H \hra \C^2$?
\end{problem}

The method of proof breaks down at the point where one would need to know that the Lie algebra 
of holomorphic vector fields vanishing on $H$ has the density property. To decide about this 
is a notoriously hard problem well known and open since decades, as is the problem about the density
property of $(\C^*)^n$ for $n>1$.

%
%

\section{Proof of Theorems \ref{th:main} and \ref{th:main2}}\label{sec:proof}

We begin by recalling some basic facts from the theory
of Stein manifolds (see e.g.\ Gunning and Rossi \cite{GunningRossi2009} or 
H{\"o}rmander \cite{Hormander1990}) and explaining the setup.

A domain $D$ in a complex manifold $Y$ is said to be {\em Runge} in $Y$ if $\{f|_D: f\in\Ocal(Y)\}$
is a dense subset of $\Ocal(D)$. If both $D$ and $Y$ are Stein,
this holds if and only if for every compact subset $K\subset D$ we have that
$\wh K_{\Ocal(D)}=\wh K_{\Ocal(Y)}$. In particular, a domain in a Stein manifold $Y$ 
which is exhausted by compact $\Ocal(Y)$-convex sets is Runge in $Y$.

A holomorphic embedding $\theta\colon X\hra Y$ of a complex manifold $X$ into a 
complex manifold $Y$ is said to be Runge if the image $Z=\theta(X)\subset Y$ is an
$\Ocal(Y)$-convex subset of $Y$, i.e., it is exhausted by compact $\Ocal(Y)$-convex subsets.
If $X$ and $Y$ are Stein manifolds, then every proper holomorphic embedding $X\hra Y$ is Runge.

Assume  that $\pi\colon E\to X$ is a holomorphic vector bundle over a Stein manifold $X$.
The total space $E$ is then also a Stein manifold. We shall write elements of $E$
in the form $e=(x,v)$ where $\pi(e)=x$, identifying $X$ with the zero section 
$\{(x,0):x\in X\}$ of $E$. For any $t\in\C$ there is a holomorphic fibre preserving map 
\begin{equation}\label{eq:psit}
	\psi_t\colon E\to E,\qquad \psi_t(x,v)=(x,tv).
\end{equation}
Clearly, $\psi_t$ is a holomorphic automorphism of $E$ for every $t\in\C^*=\C\setminus\{0\}$. 

A subset $Z\subset E$ is called {\em radial} if $\psi_t(Z)\subset Z$ holds for every $t\in[0,1]$. 

The following lemma provides the induction step in the proof of Theorem \ref{th:main}. 

%
%
\begin{lemma}\label{lem:main}
Assume that $X$ is a Stein manifold,  $E\to X$ is a holomorphic vector bundle,
$K\subset L$ are compact radial $\Ocal(E)$-convex subsets of $E$, $\Omega\subset E$
is an open set containing $X\cup K$, $Y$ is a Stein manifold with the density property
such that $\dim Y=\dim E$, and $\theta:\Omega\hra Y$ is a holomorphic embedding such that 
$\theta|_X\colon X\hra Y$ is a Runge embedding and $\theta(K)$ is $\Ocal(Y)$-convex.
Then, $\theta$ can be approximated as closely as desired uniformly on $K$ 
by a holomorphic embedding $\tilde \theta\colon \wt \Omega\hra Y$ of a domain $\wt \Omega\subset E$
with $X\cup L\subset \wt \Omega$ such that $\tilde \theta|_X\colon X\hra Y$ is a Runge embedding
and the sets $\tilde \theta(L)$ and $\tilde\theta(K)$ are $\Ocal(Y)$-convex.

If, in addition to the hypotheses above, $Y=\C^n$ with $n\ge \dim X+2$ and $\theta(X)\subset \C^n$ is a closed 
algebraic submanifold of $\C^n$, then the approximating embedding $\tilde\theta:\wt\Omega \hra \C^n$  can be chosen
to agree with $\theta$ on $X$.
\end{lemma}

The conditions in the lemma imply that $E\to X$ is the normal bundle of the embedding 
$\theta|_X\colon X\hra Y$. 

\begin{proof}
We identify $X$ with the zero section of $E$. 

Choose a compact $\Ocal(X)$-convex subset $X_0\subset X$ such that $\pi(L)\subset X_0$.
Since the embedding $\theta|_X \colon X\hra Y$ is Runge, the image
$Y_0=\theta(X_0) \subset \theta(X)$ is $\Ocal(Y)$-convex. 
Pick a compact $\Ocal(Y)$-convex neighborhood $N\subset \theta(\Omega)$ of $Y_0$
(such exists since an $\Ocal(Y)$-convex set has a basis of compact $\Ocal(Y)$-convex neighborhoods).
Thus, $N=\theta(N_0)$ for a compact set $N_0\subset \Omega$ with $X_0\subset \mathring N_0$. 

Let $\psi_t$ be defined by \eqref{eq:psit}. Since $\pi(L)\subset X_0$ and $N_0$ is a neighborhood of $X_0$ in $E$,
we can choose $\epsilon>0$ small enough such that $\psi_\epsilon(L)\subset N_0$. 
Since $L$ is $\Ocal(E)$-convex and $\psi_\epsilon\in\Aut(E)$, the set $\psi_\epsilon(L)$ is $\Ocal(E)$-convex, 
and hence afortiori $\Ocal(N_0)$-convex. Since $\theta\colon \Omega\to \theta(\Omega)$ is a biholomorphism,
it follows that the set $\theta(\psi_\epsilon(L))$ is $\Ocal(N)$-convex, 
and hence also $\Ocal(Y)$-convex (since $N$ is $\Ocal(Y)$-convex).

Consider the isotopy of injective holomorphic maps $\sigma_t$ for $t\in [\epsilon,1]$, 
defined on an open neighborhood of $\theta(K)$ in $Y$ by the condition
\begin{equation}\label{eq:sigma-t}
	\theta\circ \psi_t =  \sigma_t\circ \theta,\qquad  t\in [\epsilon,1].
\end{equation}
Note that the following hold:
\begin{itemize}
\item[\rm (a)] $\sigma_1$ is the identity map, and
\vspace{1mm}
\item[\rm (b)]  for every $t\in  [\epsilon,1]$ the compact set $\sigma_t(\theta(K)) \subset Y$ is 
$\Ocal(Y)$-convex.
\end{itemize}
Condition (b) holds because $\psi_t(K)\subset K$ is clearly $\Ocal(E)$-convex, so the set
$\sigma_t(\theta(K)) = \theta(\psi_t(K))$ is $\Ocal(\theta(K))$-convex and hence $\Ocal(Y)$-convex
(since $\theta(K)$ is $\Ocal(Y)$-convex).

By the Anders{\'e}n-Lempert-Forstneri{\v c}-Rosay-Varolin theorem
(see \cite[Theorem 2.1]{ForstnericRosay1993} for $Y=\C^n$ and 
\cite[Theorem 4.10.5]{Forstneric2017E} for the general case), the map $\sigma_\epsilon$ can
be approximated uniformly on a neighborhood of $\theta(K)$ by holomorphic automorphisms
$\phi\in \Aut(Y)$. 

Since $\psi_\epsilon(L\cup X)= \psi_\epsilon(L)\cup X \subset \Omega$ 
by the choice of $\epsilon$, there is an open neighborhood $\wt \Omega\subset E$ of $L\cup X$ 
such that $\psi_\epsilon(\wt \Omega)\subset \Omega$. We claim that the holomorphic embedding
\begin{equation}\label{eq:tildetheta}
	\tilde\theta := \phi^{-1}\circ \theta \circ\psi_\epsilon: \wt \Omega \longhookrightarrow  Y
\end{equation}
satisfies  the lemma provided that $\phi$ is chosen close enough to $\sigma_\epsilon$
on a fixed neighborhood of $\theta(K)$. 
Indeed, since the sets $\theta(\psi_\epsilon(L))$ and $\theta(\psi_\epsilon(K))$
are $\Ocal(Y)$-convex and $\phi$ is an automorphism of $Y$, the sets $\tilde\theta(L)$ and $\tilde\theta(K)$ 
are also $\Ocal(Y)$-convex. Furthermore, 
$\tilde\theta|_X = \phi^{-1}\circ \theta|_X\colon X\hra Y$ is a Runge embedding since $\theta|_X$ is. 
Finally, on the set $K$ we have  in view of \eqref{eq:sigma-t} that
\[
	\tilde\theta = \phi^{-1}\circ \theta \circ\psi_\epsilon = \phi^{-1}\circ  \sigma_\epsilon \circ \theta.
\] 
Since the map $\phi^{-1}\circ  \sigma_\epsilon$ is close to the identity on $\theta(K)$ by the choice of $\phi$, it follows 
that $\tilde \theta$ is close to $\theta$ on $K$. This proves the first part of the lemma.

Assume now that $Y=\C^n$ and that $\theta(X)=A\subset \C^n$ is a closed algebraic submanifold, 
where $n\ge \dim A+2$. By Kaliman and Kutzschebauch \cite[Theorem 6]{KalimanKutzschebauch2008IM}, 
the Lie algebra of algebraic vector fields on $\C^n$ vanishing on $A$  enjoys the algebraic density property. 
(This means that every algebraic vector field on $\C^n$ vanishing on $A$ can be expressed by 
sums and commutators of complete algebraic vector fields on $\C^n$
vanishing on $A$. Indeed, one may use shear vector fields vanishing on $A$.)
This implies (see \cite[Proposition 4.10.4]{Forstneric2017E}) 
that the flow of any algebraic vector field vanishing on $A$ can be approximated on each compact 
polynomially convex subset by holomorphic automorphisms of $\C^n$  fixing $A$ pointwise.

Note that, up to a change of the $t$-parameter, the isotopy $\psi_t$ \eqref{eq:psit}
is the flow of a holomorphic vector field $\Vcal$ on $E$, tangent to the fibres of the projection
$E\to X$ and vanishing on the zero section $X\subset E$.
Hence, the isotopy $\sigma_t$ defined by \eqref{eq:sigma-t} is also the flow of a holomorphic 
vector field $\Wcal$ on a neighborhood of $\theta(K)$ in $\C^n$ that vanishes on the algebraic submanifold 
$A=\theta(X)\subset \C^n$. (Indeed, $\Wcal=\theta_* (\Vcal)$ is the push-forward of $\Vcal$ by the embedding $\theta$.)
By Serre's Theorems A and B \cite{Serre1955AM} we can 
approximate $\Wcal$ as closely as desired on a neighborhood of the compact polynomially convex set 
$\theta(K)$ by an algebraic vector field vanishing on $A$. By what has been said above, 
this shows that the map $\sigma_\epsilon$ can be approximated uniformly on
a neighborhood of $\theta(K)$ by holomorphic automorphisms $\phi\in \Aut(\C^n)$
satisfying $\phi(z)=z$ for all $z\in\theta(X)$. The proof is now completed just as before.
In particular, we see from \eqref{eq:tildetheta} that the embedding $\tilde\theta$ agrees
with $\theta$ on $X$. 
\end{proof}

\smallskip

\noindent {\em Proof of Theorem \ref{th:main}.}
Pick an exhaustion $K_1\subset K_2\subset \cdots \subset \bigcup_{j=1}^\infty K_j=E$
by compact radial $\Ocal(E)$-convex sets. In fact, we may choose each $K_j$ of the form
\[
	K_j=\{(x,v)\in E:\varphi_1(x)\leq j\ \ \mbox{and}\  |v|_{\varphi_2}\leq j\},
\]
where $\varphi_1$ is a strongly plurisubharmonic exhaustion function on $X$ and 
$\varphi_2$ is a suitably chosen hermitian metric on $E$.
Let $\theta\colon X\hra Y$ be a holomorphic Runge embedding. 
By a theorem of Docquier and Grauert (see \cite[Theorem 3.3.3]{Forstneric2017E}) there 
is a neighborhood $\Omega_0\subset E$ of the zero section $X\subset E$ such that
$\theta$ extends to a holomorphic embedding $\theta_0:\Omega_0\hra Y$.
Set $K_0=\varnothing$. By applying Lemma \ref{lem:main} inductively, we find a sequence of open neighborhoods 
$\Omega_j\subset E$ of $K_j\cup X$ and holomorphic embeddings $\theta_j\colon \Omega_j\hra Y$ 
satisfying the following conditions for every $j\in \N$:
\begin{itemize}
\item[\rm (a)] the compact sets $\theta_j(K_j)$ and $\theta_j(K_{j-1})$ are $\Ocal(Y)$-convex, 
\item[\rm (b)] the embedding $\theta_j|_X\colon X\hra Y$ is Runge, and 
\item[\rm (c)] $\theta_j$ approximates $\theta_{j-1}$ as closely as desired on $K_{j-1}$. 
\end{itemize}
If the approximations are close enough,  the sequence $\theta_j$ converges uniformly on compacts in $E$
to a holomorphic embedding $\tilde \theta \colon E\hra Y$. 
Since $\Ocal(Y)$-convexity of a compact set in a Stein manifold $Y$ is a stable property for compact strongly 
pseudoconvex domains \cite{Forstneric1986}, and every compact $\Ocal(Y)$-convex set can be approximated
from the outside by such domains, it follows that the image of each $K_j$ remains $\Ocal(Y)$-convex in the limit 
provided that all approximations were close enough. Hence, $\tilde \theta(E)$ is a Runge domain in $Y$.  
\qed

\smallskip

\noindent {\em Proof of Theorem \ref{th:main2}.}
We follow the proof of  Theorem \ref{th:main}. By the second part of Lemma \ref{lem:main},
the sequence of embeddings $\theta_j\colon \Omega_j\hra \C^n$ can now be chosen such that,
in addition to the above, we have that $\theta_j|_X=\theta_{j-1}|_X$ holds for all $j\in\N$.
This ensures that the limit embedding $\tilde \theta=\lim_{j\to\infty}\theta_j \colon E\hra\C^n$ 
also satisfies $\tilde\theta|_X=\theta|_X$.
\qed


\subsection*{Acknowledgements}
F.\ Forstneri{\v c} is supported by the research program P1-0291 and grant J1-7256 from ARRS, Republic of Slovenia. 
E. F.\ Wold is supported by the RCN grant 240569, Norway.  


{\bibliographystyle{abbrv} \bibliography{references}}


\vspace*{5mm}
\noindent Franc Forstneri\v c

\noindent Faculty of Mathematics and Physics, University of Ljubljana, Jadranska 19, SI--1000 Ljubljana, Slovenia

\noindent Institute of Mathematics, Physics and Mechanics, Jadranska 19, SI--1000 Ljubljana, Slovenia

\noindent e-mail: {\tt franc.forstneric@fmf.uni-lj.si}

\vspace*{5mm}

\noindent Erlend F.\ Wold 

\noindent Department of Mathematics, University of Oslo, Postboks 1053 Blindern, 
NO-0316 Oslo, Norway

\noindent e-mail:  {\tt erlendfw@math.uio.no}

\end{document}